# Dimensionality reduction of discrete-time dynamical systems

Chengyi Tu

## Abstract

One of the outstanding problems in complexity science and dynamical system theory is understanding the dynamic behavior of high-dimensional networked systems and their susceptibility to transitions to undesired states. Because of varied interactions, large number of parameters and different initial conditions, the study is extremely difficult and existing methods can be applied only to continuous-time systems. Here we propose an analytical framework for collapsing N-dimensional discrete-time systems into a S+1-dimensional manifold as a function of S effective parameters with S << N. Specifically, we provide a quantitative prediction of the quality of the low-dimensional collapse. We test our framework on a variety of real-world complex systems showing its good performance and correctly identify the regions in the parameter space corresponding to the system's transitions. Our work offers an analytical tool to reduce dimensionality of discrete-time networked systems that can be applied to a broader set of systems and dynamics.

## Introduction

The study of complex dynamic systems is rapidly attracting interest within the multidisciplinary science community, with biology, ecology, environmental science and socioeconomics being some of the many areas of investigation[1-7]. To quantitatively investigate and understand their dynamic behavior, we often need to analyze their stationary state(s), i.e., equilibrium points. In general, a complex dynamic system may have multiple attractors of different types, depending on the interaction structures, parameter values and initial conditions. In high-dimensional dynamic systems, it is often impossible to characterize the extent of the domain of attraction of its stable equilibria and how the boundaries of such a domain change. Therefore, a crucial question in complexity science and dynamic system theory is identifying the factors that would prevent the state shift from desired to undesired state. In high-dimensional complex systems, it is not possible to characterize the near equilibrium phase space as a function of the (many) parameters of the system.

Although great effort has been devoted to understanding the dynamic behavior of complex systems using one-dimensional methods[8] and critical slowing down theory[9-12], the study of the factors underlying the collapse of the high-dimensional dynamic system, especially for the networked system, remains an outstanding problem. Recently, Gao *et al*.

developed an analytical tool with which it is possible to identify the natural control and state parameters of a high-dimensional networked system through mean-field approaches[13]. Tu et al. found a new condition poses effective limitations to their framework[14] and lent it to the study of the collapse or functioning of any networked systems, accounting for the full heterogeneity in node-specific self- and coupling-dynamics[15]. Laurence et al. developed a polynomial approximation to reduce complex networks based on spectral graph theory[16,17].

However, all these methods can be applied only to the continuous-time dynamic system. In real-world applications, natural and engineered complex systems are usually in discrete-time format. For example, in ecological community abundances of each species are observed in discrete time[18]; in epidemic spread number of infected and cured people are counted in discrete time[19]. Therefore, a more general framework to explore the functioning or the dimensionality reduction of discrete-time networked dynamic systems is needed to fill the gap existing between theory and real-world problems.

Here we develop a general analytical framework that can be used to reduce the dimensionality of the "order" parameter space as a function of a set of effective "control" parameters, defined as those parameters that drive the functioning (associated with specific system states) of any discrete-time networked system. Specifically, we provide quantitative predictions of the quality of the low-dimensional collapse as a function of the properties of self- and coupling-dynamics as well as interaction networks using results from random matrix theory. The framework allows us to investigate the possible occurrence of transitions (broadly defined) from a functioning stable state to an undesired one where the networked system collapses or stops functioning in the desired way.

# Results

## Dimensionality reduction of discrete-time networked dynamic systems

Consider a discrete-time networked dynamics system consisting of $N$ nodes whose states $\mathbf{x} = (x_1, \ldots, x_N)^T$ follow the dynamic equation

$$x_i[t+1] = F_i(x_i[t]) + \sum_{j}^{N} A_{ij} G_i(x_i[t], x_j[t]) \quad (1)$$

where $F_i(x_i[t])$ are the "local" dynamics at node $i$ (or "self-dynamics") and $G_i(x_i[t], x_j[t])$ are the dynamics expressing the coupling of node $i$ with its neighbor $j$ (or "coupling-dynamics"), according to the adjacency matrix $\mathbf{A} \in R^{N \times N}$, representing the interaction network of the system, with $A_{ij}$ capturing the interaction $i \leftarrow j$. Although the dynamic behavior, for example resilience, of continuous-time systems is investigated[13,15], a framework to investigate the corresponding discrete-time systems is still missing. To formulate a general framework for the analysis of the dynamic behaviors for a discrete-time version, we first define a mean-field operator $\mathcal{L}(\mathbf{x}) = \frac{1}{N}\sum_{j=1}^{N} s_j^{out} x_j / \frac{1}{N}\sum_{j=1}^{N} s_j^{out} = \frac{\langle \mathbf{s}^{out} \cdot \mathbf{x} \rangle}{\langle \mathbf{s}^{out} \rangle}$ where $\mathbf{s}^{out} = (s_1^{out},\ldots,s_N^{out})$ is the vector of the out-degree of the matrix $\mathbf{A}$; then, we characterize the effective state of the networked system using the weighted average node state $x_{eff} = \mathcal{L}(\mathbf{x})$. If the network's degree correlation is low and node activities are uniform, Eq. (1) can be reduced to

$$I(d_1,\ldots,d_S, x_{eff}) = x_{eff}[t+1] \approx \sum_{s=1}^{S} d_s * x_{eff}^{s-1}[t] \quad (2)$$

where $S = \max(m,n)$, $d_s = \begin{cases} B_{eff}^s + A_{eff} * C_{eff}^s, s \in [1, \min(m,n)] \\ A_{eff} C_{eff}^s, s \in [m+1, n], m < n \\ B_{eff}^s, s \in [n+1, m], n < m \end{cases}$ are effective parameters where $m$ and $n$ are orders of polynomial functions of Chebyshev approximation of self-dynamics $F_i(x_i[t])$ and coupling-dynamics $G_i(x_i[t], x_j[t])$ respectively; $A_{eff} = \mathcal{L}(\mathbf{s}^{in})$, $B_{eff}^s = \mathcal{L}(B^s)$, and $C_{eff}^s = \mathcal{L}(C^s)$ (see Methods). $B^k = (b_{1,k},\ldots,b_{N,k})^T$ is the column of the $k$-th term of the $m$-order Chebyshev polynomials[20,21] approximating the self-dynamics $F_i(x_i[t])$, and $C^l = (c_{1,l},\ldots,c_{N,l})^T$ is the column of the $l$-th term of the $n$-order Chebyshev polynomials approximating the coupling-dynamics $G_i(x_i[t], x_j[t])$. Therefore, we map Eq. (1), a high-dimensional networked system, into Eq. (2), a low-dimensional effective system with $\max(n,m)$ effective parameters and the state variable $x_{eff}$. We can study the dynamic behavior of the high-dimensional networked system through the behavior of $x_{eff}$ at steady state(s). In particular, the conditions for stability of steady state(s) $x_{eff}^*$ can be associated with a region expressed by the equation set

$$\begin{cases} I(d_1,\ldots,d_S, x_{eff}^*) = 0 \\ \left\| Re\left[\frac{\partial I}{\partial x_{eff}}\big|_{x_{eff}=x_{eff}^*}\right] \right\| < 1 \end{cases} \quad (3)$$

In other words, for the low-dimensional system given by Eq. (2) we can calculate analytically the resilience function $x(d_1,\ldots,d_S)$ —uniquely determined by $I(d_1,\ldots,d_S, x_{eff})$ —which represents the possible states of the system as a function of the effective parameters $d_1,\ldots,d_S$. Therefore, in order to study the stability or the existence of critical transitions in the high-dimensional system given by Eq. (1) one has to simply calculate $d_1,\ldots,d_S$ from the high-dimensional system and analyze the corresponding resilience function $x(d_1,\ldots,d_S)$ corresponding to the low-dimensional system given by Eq. (2).

## Error definition

For estimating an error of the proposed approximation, we can define projection distance from the point $(d_1,\ldots,d_S,x_{eff})$ obtained from Eq. (1) and the stationary solution of the low-dimensional resilience function $x_{eff}(d_1,\ldots,d_S)$ obtained from Eq. (3) as the error of the proposed approximation (see Fig. 1)

$$err = \left|x_{eff} - x_{eff}(d_1,\ldots,d_S)\right| \quad (4)$$

If this error distance is small, the point of numerical simulation is near the surface of the low-dimensional resilience function, which means that this framework works well.

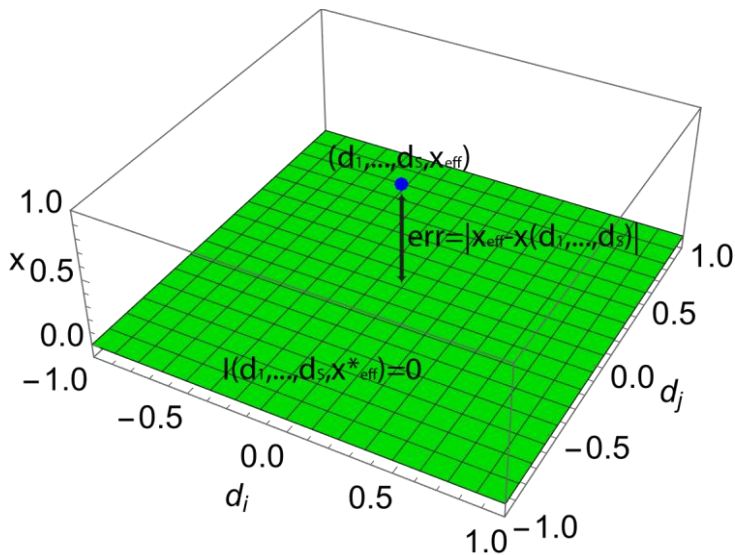

Figure 1. Quantifying the goodness of the dimensionality reduction. The blue point indicates $(d_1,\ldots,d_S,x_{eff})$ corresponding to the stationary state variables of Eq. (1). The green surface is the analytical stationary solution of the low-dimensional effective function Eq. (3). The vertical distance between the point and the curve represents the error of dimensionality reduction.

## Dynamic behavior of discrete-time GLV model with random networks

In order to test the validity of this framework, we consider a discrete-time generalized Lotka-Volterra (GLV) dynamics[22-]

[24] given by $\frac{x_i[t+1]-x_i[t]}{T} = \alpha_i x_i[t] + \sum_{j=1}^{S} A_{ij} x_i[t] x_j[t]$ where $S$ is the number of species in a community, $x_i \geq 0$ is the population size of species $i$, $\alpha_i$ is the intrinsic growth rate of species $i$ and $A_{ij}$ is the interaction between species $i$ and $j$. It is a set of first-order, non-linear, differential equations frequently used to describe the population dynamics of interacting species in community ecology. Species interaction networks can be used to model mutualistic (++), commensal (+0), competitive (--), amenalistic (-0) and predator-prey (+-) relationships among an arbitrary number of species. For simplification, we assume that $T=1$ and the above equations reads[25,26]

$$x_i[t+1] = (1+\alpha_i) x_i[t] + \sum_{j=1}^{S} A_{ij} x_i[t] x_j[t] \quad (5)$$

The vectorized form of above equation is $\mathbf{x}[t+1] = \mathbf{x}[t] \circ (1 + \boldsymbol{\alpha} + \mathbf{A} \cdot \mathbf{x}[t])$. The advantage of using GLV dynamics is that we have an analytical solution for the stationary state $\mathbf{x}^* = -\mathbf{A}^{-1}\boldsymbol{\alpha}$ or $\mathbf{x}^* = 0$ as a function of the growth rate vector $\boldsymbol{\alpha}$ and interaction network $\mathbf{A}$. The low-dimensional effective equation is $x_{eff}[t+1] = d_2 x_{eff}[t] + d_3 x_{eff}^2[t]$ where $d_2 = B_{eff}^2 = 1 + \alpha_{eff} = 1 + \mathcal{L}(\boldsymbol{\alpha})$ and $d_3 = A_{eff} = \mathcal{L}(\mathbf{s}^{in})$ and its corresponding low-dimensional resilience function is $x_{eff}(d_2, d_3) = \frac{1-d_2}{d_3}$ when $1 < d_2 < 3, d_3 < 0$ and $x_{eff}(d_2, d_3) = 0$ when $-1 < d_2 < 1$.

If all pairs of off-diagonal elements $A_{ij}$ with mean $\mu_A$ and standard deviation $\sigma_A$ do not have any correlation so that $\rho_A = 0$ and diagonal elements $D_i$ are drawn from a univariate distribution with mean $\mu_D$ and standard deviation $\sigma_D$, then the effective parameter of the adjacency matrix $\mathbf{A}$ is $A_{eff} = \frac{(\mu_D^2 + \sigma_D^2) + (S-1)\left[-2\mu_A\mu_D + (S-1)\mu_A^2\right]}{-\mu_D + (S-1)\mu_A}$ (see Methods). Therefore, $d_3 = A_{eff}$ approximates to $\frac{(\mu_D^2 + \sigma_D^2) + (S-1)\left(-2\mu_A\mu_D + (S-1)\mu_A^2\right)}{-\mu_D + (S-1)\mu_A}$. Further, if all elements of the vector $\boldsymbol{\alpha}$ are drawn from a given distribution with mean $\mu_\alpha$ and standard deviation $\sigma_\alpha$, then $d_2 = 1 + \alpha_{eff}$ approximates to $1 + \mu_\alpha$.

Assuming that $\mathbf{x}$ touches stationary state, the effective state is $x_{eff} = \frac{-\mu_\alpha}{-\mu_D + (S-1)\mu_A}$ for $\mathbf{x}^* = -\mathbf{A}^{-1}\boldsymbol{\alpha}$ or $x_{eff} = 0$ for $\mathbf{x}^* = 0$. Therefore, the error itself becomes a random variable whose probability distribution is inherited from the distributions of the random vector $\boldsymbol{\alpha}$ and random matrix $\mathbf{A}$. We get the error analytically $err = \left| \frac{\mu_\alpha \sigma_D^2}{\left((S-1)\mu_A - \mu_D\right)^2 \left(\left((S-1)\mu_A - \mu_D\right) + \sigma_D^2\right)} \right|$ for nonzero solution $x_{eff} = \frac{-\mu_\alpha}{-\mu_D + (S-1)\mu_A}$ and

$$err = \left| \frac{\mu_\alpha ((S-1)\mu_A - \mu_D)}{((S-1)\mu_A - \mu_D)^2 + \sigma_D^2} \right| \text{ for zero solution } x_{eff} = 0.$$

We then perform numerical simulations to investigate how changes in the parameters of Eq. (5) affect species abundances **x** for testing whether this framework works well. We set $\mu_\alpha = 1, \sigma_\alpha = |\mu_\alpha/3|$, $\mu_X = -0.04, \sigma_X = |\mu_X/3|$, $\mu_D = 1, \sigma_D = |\mu_D/3|$, $C = 0.5$, $S = 50$ and use two different initial conditions for $\mathbf{x}(t=0)$: a low initial population (i.e., the elements of **x** are randomly drawn from a uniform distribution between 0 and 0.1), and a high initial population (i.e., the elements of **x** are randomly drawn from a uniform between 0.9 and 1). The results can be observed in Fig. 2. The effective equation shows that for a given set of parameters the dynamics have only one stable equilibrium. We use our framework to identify the changes in the parameter space that are associated with a collapse from a state of species coexistence ($x_{eff} > 0$) to the state where all species go extinct ($x_{eff} = 0$). We find that the results from the numerical simulations collapse onto the manifold $x_{eff}(d_2, d_3)$ (Fig. 2 a) with relatively small errors. To highlight the performance, in Fig. 2 we compare the errors from analytical predictions and numerical simulations. We note again that $d_2(\alpha_{eff})$ and $d_3(A_{eff})$ depend on the specific realization (of matrix **A** and vector **α**) over the different realizations. In particular, we find that $d_2 = 1$, i.e., $\alpha_{eff} = 0$, is the critical value for such a transition. Additionally, we find $err \to 0$ as $S \to +\infty$.

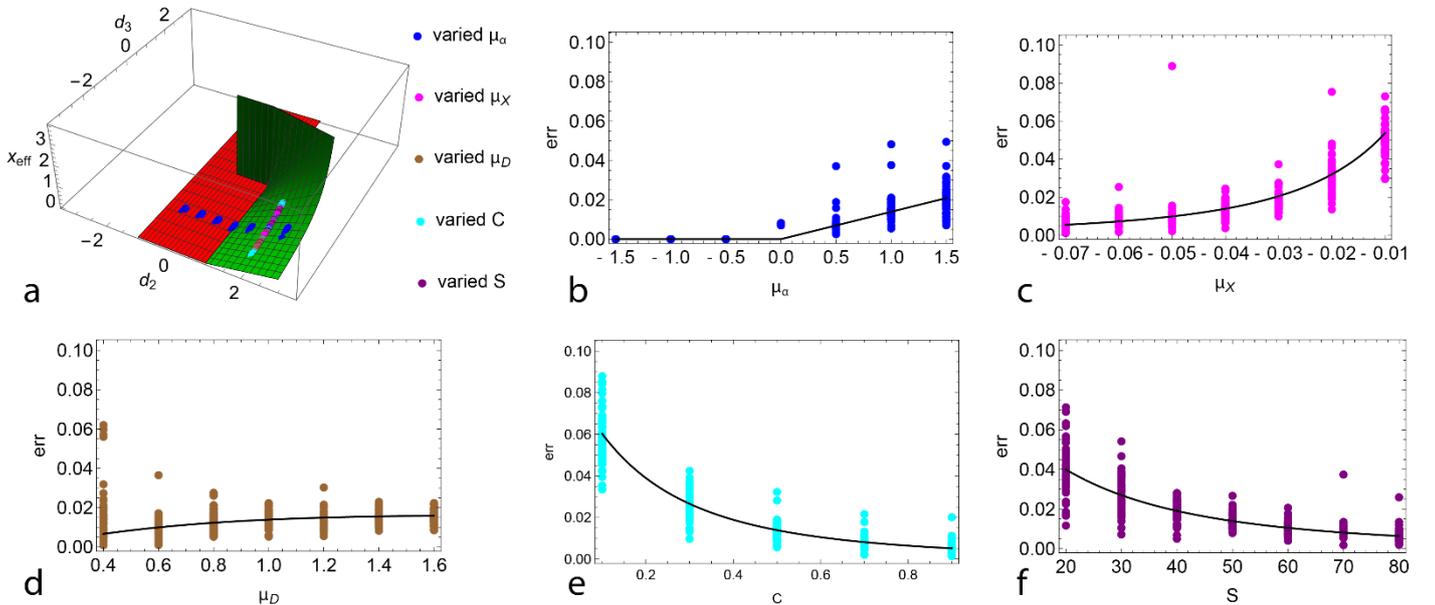

Figure 2. (a) Effective state of the system in three-dimensional space composed of the state variable $x_{eff}$ and effective parameters $d_2, d_3$ with varied $\mu_\alpha$, $\mu_X$, $\mu_D$, $C$ and $S$. Each colored surface represents one stable state in the manifold. The points representing the steady states in the high-dimensional model as a function of changes in the parameters of the dynamics collapse onto the manifold. Error as a function of changes in (b) $\mu_\alpha$, (c) $\mu_X$, (d) $\mu_D$, (e) $C$

and (f) $S$. In each panel, the solid black line represents the error of the analytical prediction. For each case, we run 50 simulations.

## Dynamic behavior of discrete-time GLV model with empirical networks

In this section, we consider empirical ecology networks where a community is composed of $S$ species comprising $S_p$ plants and $S_a$ animals such as insects serving as pollinators with $S = S_p + S_a$ [27-31]. $x_i^p$ and $x_j^a$ denote the abundances of the $i$-th plant species and the $j$-th animal pollinator species respectively and $\mathbf{x} = \{x_1^p, x_2^p, \ldots, x_{S_p}^p, x_{S_p+1}^a, \ldots, x_S^a\}$ is a vector expressing the population size of each of the species in the community. Species-specific intrinsic growth rates are the elements of the $S$-dimensional vector $\boldsymbol{\alpha}$. The species' interaction matrix $\mathbf{A}$ is composed of four blocks, two of which describe the direct competitive interactions among plants ($\Omega_{pp}$) and insects ($\Omega_{aa}$), respectively, and the other two define the mutualistic interactions between insects and plants ($\Gamma_{pa}$) and vice-versa ($\Gamma_{ap}$). Therefore, the interaction matrix $\mathbf{A}$ has the following structure $\begin{bmatrix} \Omega_{pp} & \Gamma_{pa} \\ \Gamma_{ap} & \Omega_{aa} \end{bmatrix}$ [28]. For each plant and animal species, we set competition coefficients $\beta$ (in the $\Omega_{pp}$ and $\Omega_{aa}$ matrices) sampled from a uniform distribution with maximum -0.001 and mean $-1/S^p$ and $-1/S^a$. In other words, the interaction matrix $\mathbf{A}$ exhibits a mixture of positive and negative signs and also correlation among elements, and where intraspecific competition coefficient is set equal to -1. $\Gamma_{pa}$ and $\Gamma_{ap}$ matrices describe how species interact mutualistically. We express the weights of this interaction matrix using a trade-off function that defines the mutualistic dependence between species $j$ and $i$ as a function of their degree: $\gamma_{ij} = \gamma y_{ij} / k_i$ where $\gamma$ is drawn from a normal distribution with mean $\mu_\gamma$ and standard deviation $\sigma_\gamma = |\mu_\gamma / 3|$, $k_i$ is the degree of species $i$, and $y_{ij} = 1$ if species $i$ and $j$ interact and zero otherwise[28]. The adjacency matrix $\mathbf{Y}$ is taken directly from empirical plant-pollinator networks from different regions of the world using the Web of Life database (http://www.web-of-life.es). We only include 134 networks that have <200 species to limit the computational cost of our analyses[5]. We set $\mu_\alpha = 1, \sigma_\alpha = |\mu_\alpha / 3|$, $\mu_\gamma = 0.4, \sigma_X = |\mu_\gamma / 3|$ and use two different initial conditions for $\mathbf{x}(t=0)$: a low initial population (i.e., the elements of $\mathbf{x}$ are randomly drawn from a uniform distribution between 0 and 0.1), and a high initial population (i.e., the elements of $\mathbf{x}$ are randomly drawn from a uniform between 0.9 and 1).

We then perform numerical simulations to investigate how changes in the parameters and networks affect species abundances $\{x_1, x_2, \ldots, x_S\}$ and corresponding errors. The panels a-b of Fig. 3 show the numerical simulation of the GLV dynamics with empirical networks when: 1) decreasing the growth rate of species, mimicking the decrease of resources or increase of unsuitable environmental conditions; 2) decreasing the strength of mutualistic interactions, mimicking the possible effect of climate change[32,33]. We find that the results from the numerical simulations collapse onto the manifold $x_{eff}(d_2, d_3)$ with relatively small errors. To highlight the approximation errors, in Fig. 3 c-d we also show the errors from numerical simulations for all empirical networks. We note again that $d_2(\alpha_{eff})$ and $d_3(A_{eff})$ depend on the specific realization (of matrix $\mathbf{A}$ and vector $\boldsymbol{\alpha}$) over the different realizations. In particular, we find that $\alpha_{eff} = 0$ is the critical value for such a transition. On the other hand, a decrease in $\mu_\gamma$ leads to a decrease in $A_{eff}$ which is associated with the extinction of some of the species and a slow decrease of both $<x>$ and $x_{eff}$. Therefore, both an increase in intrinsic growth rates ($\mu_\alpha$) and mutualistic strengths ($\mu_\gamma$) are beneficial for species' coexistence.

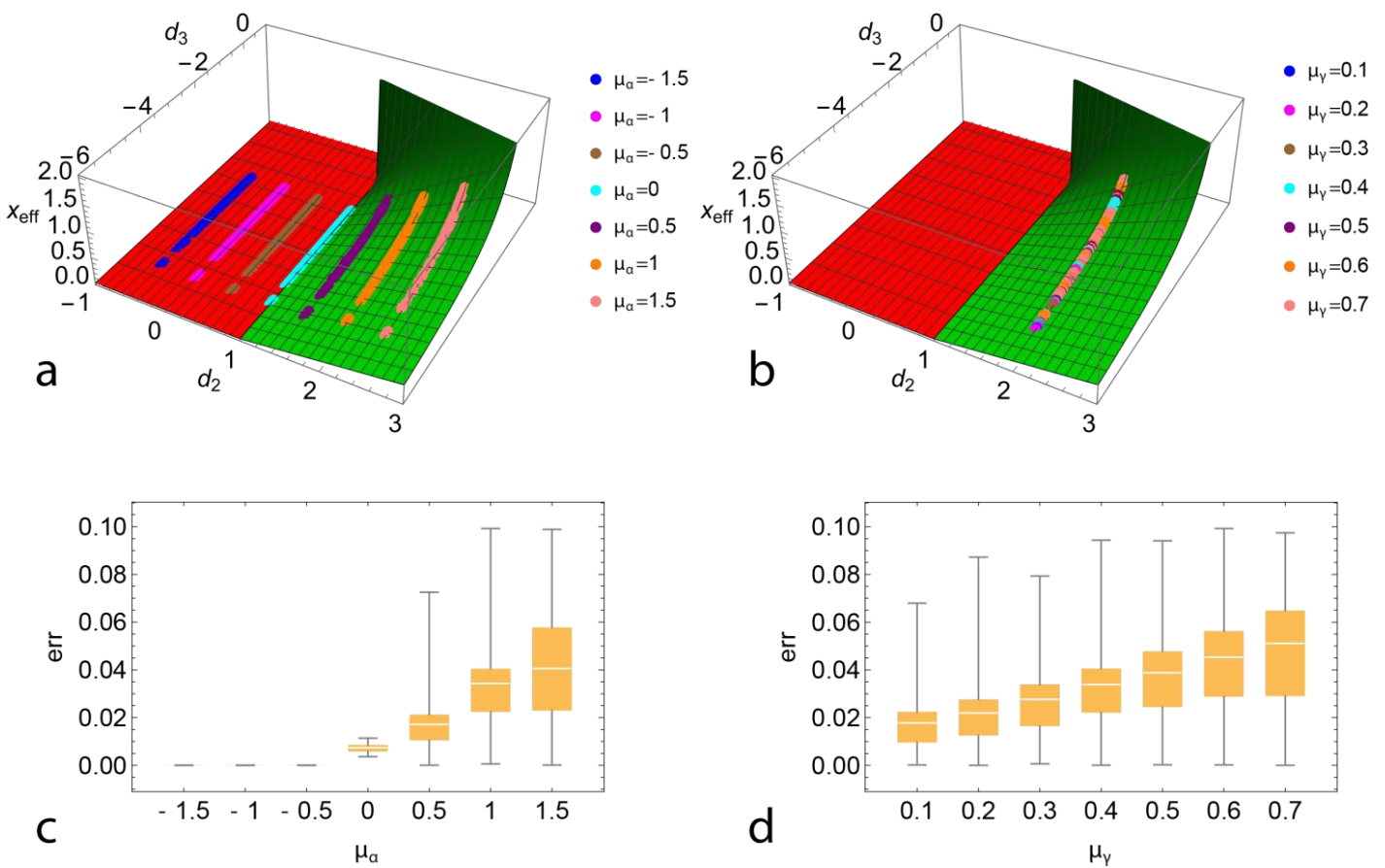

Figure 3. Effective state of the system in three-dimensional space composed of the state variable $x_{eff}$ and effective parameters $d_2, d_3$ with varied (a) $\mu_\alpha$ and (b) $\mu_\gamma$. Each colored surface represents one stable state in the manifold. The points representing the steady states in the high-dimensional model as a function of changes in the parameters of the dynamics collapse onto the manifold. Error as a function of changes in (c) $\mu_\alpha$ and (d) $\mu_\gamma$ corresponding to the panel (a)

## Dynamic behavior of discrete-time SIS model with social networks

The framework also works beyond GLV dynamics. We consider a system of $S$ individuals with some of them infected by a virus or other transmissible disease that spreads through the system as individuals interact. In the dynamics, it is often important to understand the effectiveness of various measures in limiting the spread of the epidemic. We consider the dynamics of a commonly-used susceptible-infected-susceptible (SIS) model[19,34], governed by $\frac{x_i[t+1]-x_i[t]}{T} = -e_i x_i + \sum_j^N A_{ij}(1-x_i[t])x_j[t]$. For simplification, we assume that $T=1$ and the above equation reads

$$x_i[t+1] = (1-e_i)x_i[t] + \sum_j^N A_{ij}(1-x_i[t])x_j[t] \quad (6)$$

where $0 \leq x_i \leq 1$ denotes the probability that node $i$ is in the infected state, $e_i$ is the recovery rate of node $i$ and $A_{ij}$ represents the infection rate of node $i$ as a result of the interaction with node $j$. The first term on the right-hand side of Eq. (6) accounts for the process of recovery, and the second term accounts for the process of infection. Considering the case of parameter values (i.e., $\mathbf{e}$ and $\mathbf{A}$) for which the steady state of the system is positive (i.e., $\mathbf{x}^*(t \to \infty) > 0$ or "epidemic active phase"), we can investigate how such a stationary state changes as a result of network structure.

The effective equation of the SIS model is $x_{eff}[t+1] = d_2 x_{eff}[t] + d_3 x_{eff}^2[t]$ where $d_2 = B_{eff}^2 + A_{eff} = 1 - e_{eff} + A_{eff} = 1 - \mathcal{L}(\mathbf{e}) + \mathcal{L}(\mathbf{s}^{in})$, $d_3 = -A_{eff} = -\mathcal{L}(\mathbf{s}^{in})$. It is easy to see from this low-dimensional equation that the dynamics have two steady states: $x_{eff}(d_2, d_3) = \frac{1-d_2}{d_3}$ when $1 < d_2 < 3, d_3 < 0$ and $x_{eff}(d_2, d_3) = 0$ when $-1 < d_2 < 1$. In other words, our framework predicts that the steady states of this high-dimensional equation (6) are two surfaces, $x_{eff}^* = \frac{1-d_2}{d_3}$ and $x_{eff}^* = 0$ in the space $(x_{eff}, d_2, d_3)$, and their stability depends on the order parameters $d_2$ and $d_3$.

We test these predictions by simulating the dynamics and exploring the effect of changes in the parameters characterizing the network structures. Many empirical social networks have Erdős–Rényi (ER), Barabasi-Albert (BA) or Small-World (SW) structure[35], so we use them in numerical simulation to model real-world social networks. ER network means each edge has a fixed probability of being present or absent, independently of the other edges. It can be used in the probabilistic method to prove the existence of graphs satisfying various properties. BA network means it has scale-free structure of the network degree distribution (i.e., how the number of the different node's connections are distributed). The origin of this distribution is quite common in social network, as it can be obtained through a preferential attachment mechanism. Scale-free means that few players (the "hubs") have extremely large neighbors, while most of the players have very few neighbors. SW network means it includes short average path lengths and high clustering so that most players can "be reached" from every other player by a small number of steps, and transitivity in the relation between players is observed. We express the recovery rates $\mathbf{e} = (e_1,\ldots,e_N)^T$ as random parameters drawn from a uniform distribution between 0 and $2\mu_e$ where $\mu_e$ is thus the mean recovery rate, and default value $\mu_e = 0.5$. We set two initial conditions for $\mathbf{x}(t=0)$: a low contagion initial state whose elements are drawn from a uniform distribution between 0 and 0.1, and a high contagion initial state whose elements are drawn from a uniform distribution between 0.9 and 1.

Fig. 4 shows the results for $x_{eff} = \frac{\langle \mathbf{s}^{out} \cdot \mathbf{x} \rangle}{\langle \mathbf{s}^{out} \rangle}$ obtained from numerical simulation of Eq. (6) as a function of the following changes in the model's parameters: 1) we change parameter of ER network, i.e., network connectivity, 2) we change parameter of BA network, i.e., added edges of each new node, and 3) we change parameter of SW network, i.e., rewiring probability of initial regular graph. We find that the results from the numerical simulations collapse onto the manifold $x_{eff}^*(d_2, d_3)$ (Fig. 4 a) with a relatively small error and initial conditions do not play any role. To better appreciate the approximation error, in Fig. 4 b-d we also show the error from numerical simulations. We highlight that $d_2(e_{eff}, A_{eff})$ and $d_3(A_{eff})$ depend on the specific realization (of matrix $\mathbf{A}$ and vector $\mathbf{e}$) over the different realizations.

We can thus see that small changes in $d_2$ may lead to the collapse of the epidemic, i.e. a transition from $x_{eff} > 0$ to $x_{eff} = 0$. Specifically, if $A_{eff} > e_{eff}$ then $d_2 > 0$ and $x_{eff} > 0$. Vice versa, if $A_{eff} < e_{eff}$ then $d_2 < 0$ and $x_{eff} = 0$. As the recovery rate $\mu_e$ increases, $d_2$ decreases, thereby leading to a decrease of $x_{eff}$ to zero. Similarly, as mutual interactions are reduced, $A_{eff}$ decreases, leading to a decrease in both $d_2, d_3$, thereby driving $x_{eff}$ towards the collapse. We highlight that the presented results are independent of the system's initial conditions and our framework works for different network structures.

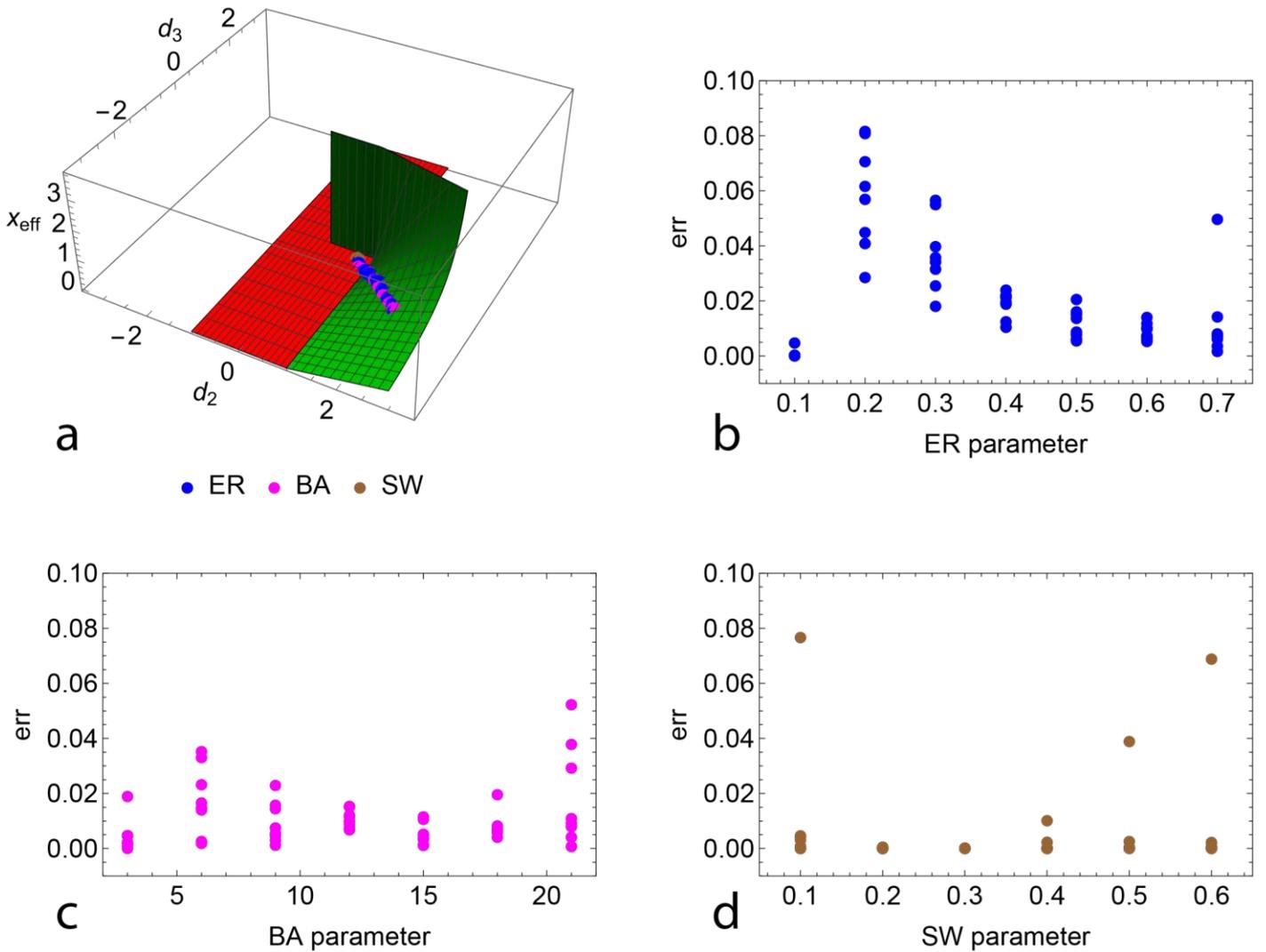

Figure 4. (a) Effective state of the system in three-dimensional space composed of the state variable $x_{eff}$ and effective parameters $d_2, d_3$. Each colored surface represents one stable state in the manifold. The points representing the steady states in the high-dimensional model as a function of changes in the parameters of the dynamics collapse onto the manifold. Error as a function of changes in (b) parameter of ER network, (c) parameter of BA network and (d) parameter of SW network. For each case, we run 50 simulations.

# Discussion

The dynamic behavior of networked systems is always difficult to investigate, although great effort has been devoted to it. The existing methods can be applied only to the continuous-time dynamical system, while the natural and engineered systems are usually in discrete-time format in real world. Therefore, a more general framework to explore the functioning or the dimensionality reduction of discrete-time networked dynamical systems is needed to fill the gap existing between theory and real-world problems. Our framework can map a high-dimensional networked system into a low-dimensional effective system in only one effective state variable and $S$ effective parameters and it works well in different dynamic models, network structures as well as node-specific self- and coupling-dynamics. We can then use these manifold dynamics to investigate the system's response to changes in the interaction structures, parameter values and initial conditions, determine the steady states of the system and the possible existence of transitions between functioning or sustainable stables including possible critical transitions in the case of bifurcating bi-stable dynamics. The numerical simulation of real-world examples presented in the Results section demonstrates that our framework works well, especially for social networks with ER, BA, or SW structure. Furthermore, we provide a quantitative prediction of the quality of the low-dimensional collapse as a function of the properties of self- and coupling-dynamics as well as interaction networks using results from random matrix theory. In summary, we show how the average dynamics of a high-dimensional networked system can be captured by a low-dimensional manifold characterizing the role of interaction network, self-dynamics and coupling-dynamics in the equilibrium states of the system and their dependence on the system's parameters. Our methodology provides (approximated) results for dimensionality reduction that are applicable to a broader set of systems and dynamics, and exploits in different contexts ranging from ecology to epidemiology and the study of critical transitions.

# Methods

## Reduction of high-dimensional dynamical system

As given by Eq. (1) of main text, the dynamics of each node depend on the node itself (given by the "self-dynamics" $F_i(x_i[t])$ ) and the interaction with its nearest neighbors (given by the interaction network and coupling dynamics $\sum_{j=1}^{N} A_{ij} G_i(x_i[t], x_j[t])$ ). Therefore, the dynamics of the average nearest-neighbor nodes represent an important contribution to the overall system's dynamics. To quantify this contribution, we define an operator $\mathcal{L}(\mathbf{x}) = \frac{1}{N}\sum_{j=1}^{N} s_j^{out} x_j / \frac{1}{N}\sum_{j=1}^{N} s_j^{out} = \frac{\langle \mathbf{s}^{out} \cdot \mathbf{x} \rangle}{\langle \mathbf{s}^{out} \rangle}$ where $\mathbf{s}^{out} = (s_1^{out}, \ldots, s_N^{out})$ is the vector of the out-degree of the interaction network $\mathbf{A}$ [13]. The operator $\mathcal{L}$ is feasible for linear time-invariant (LTI) function, i.e., $\mathcal{L}(a\mathbf{x}+b\mathbf{y}) = a\mathcal{L}(\mathbf{x}) + b\mathcal{L}(\mathbf{y})$, and Hadamard product, i.e., $\mathcal{L}(\mathbf{x} \circ \mathbf{y}) \approx \mathcal{L}(\mathbf{x})\mathcal{L}(\mathbf{y})$ [15].

If the degree correlation of network $\mathbf{A}$ is weak (the neighborhood of node $i$ is similar to the neighborhoods of all other nodes), then $\sum_{j}^{N} A_{i,j} G_i(x_i[t], x_j[t]) \approx s_i^{in} \mathcal{L}(G_i(x_i[t], \mathbf{x}[t]))$. Furthermore, if $G_i(x_i[t], x_j[t])$ is linear in $x_j$ or the standard deviation in the elements of the vector $\mathbf{x}$ is small, then $\mathcal{L}(G_i(x_i[t], \mathbf{x}[t])) \approx G_i(x_i[t], \mathcal{L}(\mathbf{x}[t])) = G_i(x_i[t], x_{eff}[t])$ where $x_{eff}[t] = \mathcal{L}(\mathbf{x}[t])$. Therefore, Eq. (1) in the main text can be written as $x_i[t+1] \approx F_i(x_i[t]) + s_i^{in} G_i(x_i[t], x_{eff}[t])$, and its vector notation is

$\mathbf{x}[t+1] = \mathbf{F}(\mathbf{x}[t]) + \mathbf{s}^{in} \circ \mathbf{G}(\mathbf{x}[t], x_{eff}[t])$ where $\mathbf{F}(\mathbf{x}[t]) = (F_1(x_1[t]), \ldots, F_N(x_N[t]))^T$ and

$\mathbf{G}(\mathbf{x}[t], x_{eff}[t]) = (G_1(x_1[t], x_{eff}[t]), \ldots, G_N(x_N[t], x_{eff}[t]))^T$.

If each $F_i(x_i[t])$ is a linear combination of $m$ subfunctions, i.e.,

$F_i(x_i[t]) = b_{i,1} f_1(x_i[t]) + b_{i,2} f_2(x_i[t]) + \cdots + b_{i,m} f_m(x_i[t])$, then according to the feasible of operator $\mathcal{L}$ for LTI and

Hadamard product, $\mathcal{L}(\mathbf{F}(\mathbf{x}[t])) \approx \sum_{k=1}^{m} \mathcal{L}(B^k) f_k(x_{eff}[t])$ where $B^k = (b_{1,k}, \ldots, b_{N,k})^T$ is the $k$-th column of matrix

$B$. Generally, we can use Chebyshev polynomials to approximate $F_i(x_i)$ to $\sum_{k=1}^{m} b_{i,k} x^{(k-1)}$, minimizing the error between them. Therefore, $F_i(x_{eff}[t]) = \sum_{k=1}^{m} b_{i,k} x_{eff}^{(k-1)}[t]$ [15]. Similarly, if each $G_i(x_i[t], x_j[t])$ is a linear combination of $n$ subfunctions, i.e., $G_i(x_i[t], x_j[t]) = c_{i,1} g_1(x_i[t], x_j[t]) + \cdots + c_{i,n} g_n(x_i[t], x_j[t])$, then

$\mathcal{L}(\mathbf{G}(\mathbf{x}[t], x_{eff}[t])) \approx \sum_{l=1}^{n} \mathcal{L}(C^l) g_l(x_{eff}[t], x_{eff}[t])$ where $C^l = (c_{1,l}, \ldots, c_{N,l})^T$ is the $l$-th column of matrix $C$.

Generally, we can use Chebyshev polynomials to approximate $G_i(x_i[t], x_j[t])$ to $\sum_{p,q=1}^{n/2} d_{p,q} x_i^{(p-1)}[t] x_j^{(q-1)}[t]$. Therefore,

$G_i(x_{eff}[t], x_{eff}[t]) = \sum_{l=1}^{n} c_{i,l} x_{eff}^{(l-1)}[t]$ where $c_{i,l}$ collects all terms $d_{p,q}$ such that $l = p + q - 1$.

We apply the operator $\mathcal{L}$ to both sides of the vector notation and obtain effective equation

$x_{eff}[t+1] \approx \sum_{k=1}^{m} B_{eff}^k x_{eff}^{(k-1)}[t] + A_{eff} \sum_{l=1}^{n} C_{eff}^l x_{eff}^{(l-1)}[t]$ where $A_{eff} = \mathcal{L}(\mathbf{s}^{in})$, $B_{eff}^k = \mathcal{L}(B^k)$ and $C_{eff}^l = \mathcal{L}(C^l)$ where

$k = 1, \ldots, m$ and $l = 1, \ldots, n$. To further decrease the number of parameters, we can rewrite it as

$I(d_1, \ldots, d_S, x_{eff}) = x_{eff}[t+1] \approx \sum_{s=1}^{S} d_s * x_{eff}^{s-1}[t]$ where $S = \max(m, n)$, $d_s = \begin{cases} B_{eff}^s + A_{eff} * C_{eff}^s, s \in [1, \min(m,n)] \\ A_{eff} C_{eff}^s, s \in [m+1, n], m < n \\ B_{eff}^s, s \in [n+1, m], n < m \end{cases}$ [15].

## Random matrix theory

In the most general case, we consider a matrix $\mathbf{A}$ where all pairs of off-diagonal elements - $A_{ij}$ and $A_{ji}$ - are drawn from a bivariate distribution with mean $\mu_A$, standard deviation $\sigma_A$ and correlation coefficient $\rho_A$ as well as diagonal elements $A_{ii} = -D_i$. The distribution from which the elements are drawn is not important, as only mean, variance and correlations are the relevant parameters[36-38]. Under this setting one can generate both directed and undirected networks, being able to tune also the interaction properties[38,39]. The following approximate equations would strictly hold only in the very large $S$: $\mu_A = \frac{1}{S(S-1)} \sum_{i \neq j} A_{ij}$, $\sigma_A^2 = \frac{1}{S(S-1)} \sum_{i \neq j} A_{ij}^2 - \mu_A^2$ and $\rho_A = \frac{\frac{1}{S(S-1)} \sum_{i \neq j} A_{ij} A_{ji} - \mu_A^2}{\sigma_A^2}$ where $S$ is the matrix size. According to the definition $A_{eff} = \frac{\sum_{ijk} A_{ik} A_{kj}}{\sum_{ij} A_{ij}}$ and $x_{eff} = \frac{\sum_{ij} A_{ij} x_j}{\sum_{ij} A_{ij}}$, the following approximations for $A_{eff}$ are obtained $A_{eff} = \frac{\sum_i (-D_i)^2 + (S-1) \left[ 2\mu_A \sum_i (-D_i) + S(S-1)\mu_A^2 + S\rho_A \sigma_A^2 \right]}{\sum_i (-D_i) + S(S-1)\mu_A}$ and $x_{eff} = \frac{\sum_{ij} A_{ij} x_j}{\sum_i (-D_i) + S(S-1)\mu_A}$ [6,14].

In the previous analysis, we assume that connectivity (the fraction of non-zero elements) of adjacency matrix $\mathbf{A}$ is always one, i.e., $C = 1$. Generalizing the results to not fully connected networks is straightforward. If the mean and standard deviation of the non-zero elements are $\mu_X$ and $\sigma_X$ respectively, then the global mean and standard deviation become $\mu_A = C\mu_X$ and $\sigma_A = \sqrt{C\sigma_X^2 + C(1-C)\mu_X^2}$.

## Data and code availability

The ready-to-use notebook codes to reproduce the results presented in the current study are available in OSF with the access code 5hmsn (https://osf.io/5hmsn/).